\documentclass[]{theclass}
\usepackage{inputenc}

\begin{document}
\begin{frontmatter}

\titledata{Even orientations and Pfaffian graphs}{}           
\authordata{Mari\`{e}n Abreu}
{Dipartimento di Matematica, Informatica ed Economia\\ Universit\`{a} degli Studi della Basilicata, Italy}{marien.abreu@unibas.it}{}
{}

\authordata{Domenico Labbate}
{Dipartimento di Matematica, Informatica ed Economia\\ Universit\`{a} degli Studi della Basilicata, Italy}{domenico.labbate@unibas.it}{}
{}

\authordata{Federico Romaniello}
{Dipartimento di Matematica, Informatica ed Economia\\ Universit\`{a} degli Studi della Basilicata, Italy}{federico.romaniello@unibas.it}{}
{}

\authordata{John Sheehan}
{ Department of Mathematical Sciences\\King's College, Scotland}
{j.sheehan@maths.abdn.ac.uk}{}


\begin{abstract}
We give a characterization of Pfaffian graphs in terms of even orientations,
extending the characterization of near bipartite non-pfaffian graphs
by Fischer and Little \cite{FL}. Our graph theoretical characterization is equivalent to the one proved by Little in \cite{L73} (cf. \cite{LR}) using linear algebra arguments.\end{abstract}

\end{frontmatter}

\section{Introduction}
All graphs considered are finite and simple (without loops or
multiple edges) unless otherwise stated. 
Most of our terminology is standard and can be found in many textbooks such as \cite{BM} and \cite{LP}.

Let $F$ be a $1$-factor of a graph $G$. Then a cycle $C$ is said to be {\em $F$-alternating} if $|E(C)|= 2|E(F) \cap E(C)|$. In particular, each $F$-alternating cycle has an even number of edges. An $F$-alternating cycle $C$ in an orientation $\vec{G}$ of $G$ is {\em evenly (oddly) oriented} if for either choice of direction of traversal around $C$, the number of edges of $C$ directed in the direction of traversal is even (odd). Since $C$ is even, this is clearly independent of the initial choice of direction around $C$.

Let $\vec{G}$ be an orientation of $G$ and $F$ be a $1$-factor of $G$. If every $F$-alternating cycle is evenly oriented then $\vec{G}$ is said to be an {\em even $F$-orientation of $G$}. On the other hand if every $F$-alternating cycle is oddly oriented then $\vec{G}$ is said to be an {\em odd $F$-orientation of $G$}.

An even subdivison of $G$ is any graph $G^{*}$ which can be obtained from $G$ by replacing edges $(u,v)$ of $G$ by paths $P(u,v)$ of odd length, such that $V(P(u,v)) \cap V(G) =\left\{u,v\right\}$.

An $F$-orientation $\vec{G}$ of a graph $G$ is {\em Pfaffian} if it is odd.
It turns out that if $\vec{G}$ is a Pfaffian $F$-orientation then  $\vec{G}$ is a Pfaffian $F^*$-orientation for all $1$-factors $F^*$ of $G$
(cf.\cite[Theorem 8.3.2 (3)]{LP}). In this case we simply say that $G$ is {\em Pfaffian}.

It is well known that every planar graph is Pfaffian and that the smallest non-Pfaffian graph is the complete bipartite graph $K_{3,3}$ (cf.\cite{CLM}). The Petersen graph is a further example of a non-Pfaffian graph (see \cite[Section 3]{ALS} for details).

The literature on Pfaffian graph is extensive and the results often profound (see \cite{T} for a complete survey).
In particular, the problem of characterizing Pfaffian bipartite graphs was posed  by P\'olya \cite{Po13}.
Little \cite{L} obtained the first such characterization in terms of a
family of forbidden subgraphs. Unfortunately, his characterization does not
give rise to a polynomial algorithm for determining whether a given
bipartite graph is Pfaffian, or for calculating the permanent of its
adjacency matrix when it is. Such a characterization was subsequently
obtained independently by McCuaig \cite{McC,McC1}, and Robertson, Seymour and
Thomas \cite{RST}. As a special case their result gives a polynomial algorithm, and hence a good characterization, for
determining when a balanced bipartite graph $G$ with adjacency matrix $A$ is {\em det-extremal} i.e. it has $|det(A)|=per(A)$.
For a structural characterization of det-extremal cubic bipartite graphs the reader may also refer to \cite{Thom86}, \cite{McC00}, \cite{McC1} and \cite{FJLS}.

The problem of characterizing Pfaffian general graphs seems much harder. Nevertheless, some very interesting connections in terms of {\em bricks} and {\em near bipartite graphs} have been found (cf. e.g. \cite{FL}, \cite{LP}, \cite{NT}, \cite{T}, \cite{VY}).



A graph $G$ is said to be {\em $1$-extendable} if each edge of $G$ is contained in at least one $1$- factor of $G$.
A subgraph $J$ of a graph $G$ is {\em central} if $G-V(J)$ has a 1-factor.

A $1$-extendable non-bipartite graph $G$ is said to be {\em near bipartite} if there exist edges $e_1$ and $e_2$ such
that $G \backslash \{e_1 , e_2\}$ is $1$-extendable and bipartite.

The Pfaffian property which holds for odd $F$-orientations
does not hold for even $F$-orientations. Indeed, the Wagner graph $W$ (cf. Section \ref{preliminaries}) is Pfaffian,
so there is an odd orientation which works for all $1$-factors. On the other hand, it has an even $F_1$-orientation and no even
$F_2$-orientation where $F_1$ and $F_2$ are chosen $1$-factors of $W$ (cf. Lemma \ref{claim2.6}).

A graph $G$ is said to be \emph{simply reducible to a graph} $H_0$ if $G$ has an odd length
cycle $C$ such that $H_0$ can be obtained from $G$ by contracting $C$, i.e. by removing an edge from the graph while simultaneously merging the two vertices that it previously joined and disregarding loops or multiple edges. More generally $G$ is said to be
\emph{reducible to a graph} $H$ if for some fixed integer $k$ there exist graphs $G_0, G_1, \ldots, G_k$
such that $G_0=G$, $G_k=H$ and for $i$, $1 \le i \le k$, $G_{i-1}$ is simply reducible to $G_i$.

\begin{figure}[h]
\begin{center}
\includegraphics[scale=0.5]{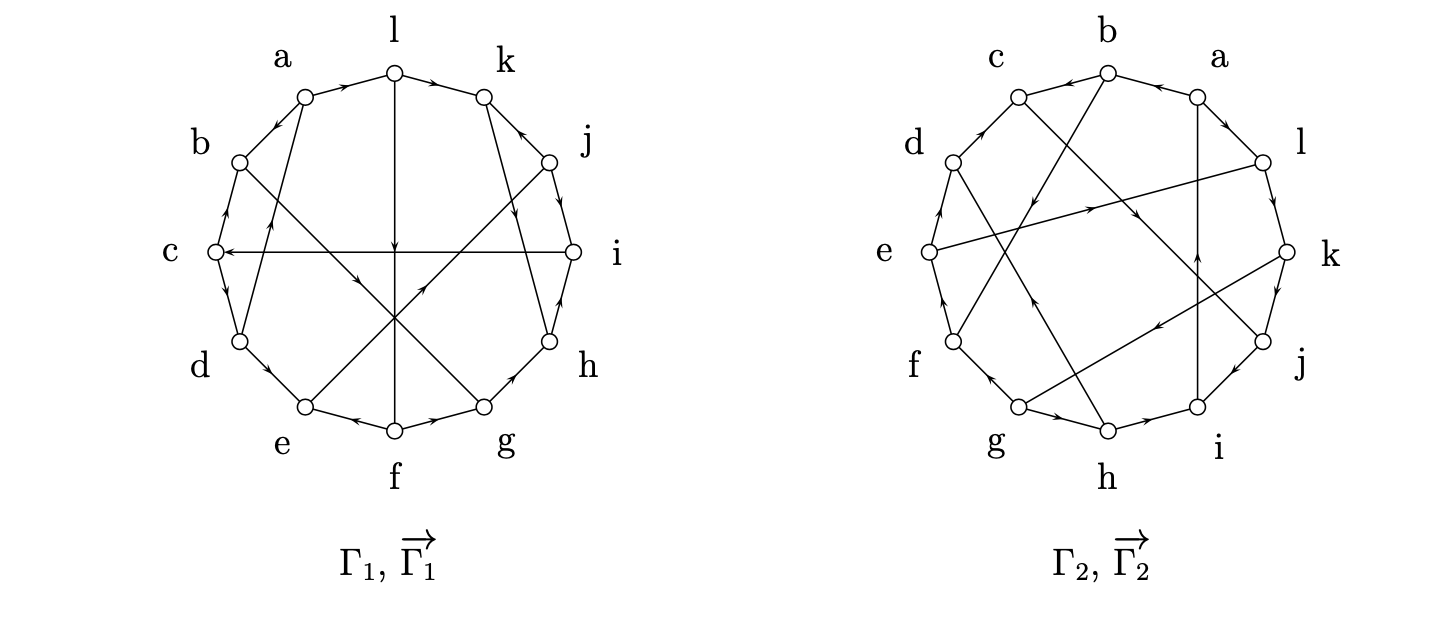}
\end{center}
  \caption{The graphs $\Gamma_1$, $\Gamma_2$ and their orientations  }\label{fig1}
\end{figure}

Fischer and Little \cite{FL} proved the following characterization of near bipartite non-Pfaffian graphs:

\begin{theorem}\cite{FL}\label{thm1_5}
A near bipartite graph $G$ is non-Pfaffian if and only if $G$ contains a central subgraph $J$ which is
reducible to an even subdivision of $K_{3,3}$, the cubeplex $\Gamma_1$ or the twinplex $\Gamma_2$ (cf. Fig. \ref{fig1} in Section \ref{bad}) \qed
\end{theorem}

In \cite{NT} this result was restated in terms of matching minors.



In this context, recently we have examined the structure of $1$-extendable graphs $G$ which have no even $F$-orientation \cite{ALS},
where $F$ is a fixed $1$-factor of $G$. We have given in \cite{ALS} a characterization in the case of graphs of connectivity at least four and of $k$-regular graphs, $k\ge 3$.
Part of this characterization is stated in Theorem \ref{lem1_6}.

In this note, as a consequence of the cited characterization of graphs with no even $F$-orientations,
we characterize non-Pfaffian graphs in terms of even orientations (cf. Theorem \ref{BadNonPfaf}),
extending the characterization of near bipartite non-Pfaffian graphs by Fischer and Little \cite{FL} cited in Theorem  \ref{thm1_5}.

Note that 
Theorem \ref{BadNonPfaf}
gives a graph theoretical proof of an equivalent formulation 
that is stated in Little and Rendl \cite{LR} and proved using linear algebra arguments in \cite{L73}.



%
%

\section{Preliminaries}\label{preliminaries}

In this section we introduce some definitions and notation useful to state and then prove our main Theorem \ref{BadNonPfaf}.



\begin{definition}\label{zerosum} {\em (Zero-sum sets)}

Let $G$ be a graph with a $1$-factor $F$.
Suppose that ${\cal A} := \{ C_1,\ldots,C_k \}$ is a set of $F$-alternating cycles such that each edge
of $G$ is contained in an even number of elements of $\cal A$.
Then ${\cal A}$ is said to be a {\em zero-sum $F$-set}.

If $k$ is even or odd we say that the zero-sum $F$-set is respectively an {\em even $F$-set} or an {\em odd $F$-set}.
\end{definition}

\begin{lemma}\label{claim2.1.1}
Let $G$ be a graph with a $1$-factor $F$ and an odd zero-sum $F$-set $\cal{C}$$:=\{C_1,\ldots,C_k\}$. Suppose that $C_1,\ldots,C_{k_1}$ are oddly $F$-oriented and $C_{k_1+1},\ldots,C_k$ are evenly $F$-oriented in an orientation $\vec{G}$ of $G$. Let $k_2:=k-k_1$ . Then, if $k_1$ is odd or $k_2$ is odd,  $G$ cannot have respectively an even $F$-orientation or an odd $F$-orientation.
\end{lemma}

\begin{proof}
Firstly suppose that $k_1$ is odd and that $G$ has an even $F$-orientation. Then there exists a set $S$ of edges such that $|E(C_i)\cap S| \equiv 1 \, (\mbox{mod} \, 2)$, $i=1,\ldots,k_1$ and $|E(C_j)\cap S| \equiv 0 \, (\mbox{mod} \, 2)$, $j=k_1+1,\ldots,k$. This follows since to change $\vec{G}$ into an even $F$-orientation we must reverse an odd number of orientations in the oddly oriented $F$-cycles and an even number of orientations in the evenly oriented $F$-cycles. Set $S:=\{e_1,\ldots,e_l\}$ and write
$$
a_{i,j} \, := \, \left \{
\begin{array}{ll}
  1 & \mbox{if} \,\, e_i \in E(C_j) \quad (j=1,\ldots,k) \\
  0 & \mbox{otherwise}
\end{array}
\right.
$$
Then, since $\cal C$ is a zero-sum $F$-set
\begin{equation}\label{eq1.1}
\sum_{j=1}^{k} a_{i,j} \,  \equiv \, 0 \, (\mbox{mod} \, 2) \,, \quad i=1,\ldots, l 
\end{equation}
and, from the definition of $S$,
\begin{equation}\label{eq1.2}
\sum_{i=1}^{l} a_{i,j} \,  \equiv \, 1 \, (\mbox{mod} \, 2) \,, \quad j=1,\ldots, k_1 
\end{equation}
\begin{equation}\label{eq1.3}
\sum_{i=1}^{l} a_{i,j} \,  \equiv \, 0 \, (\mbox{mod} \, 2) \,, \quad j=k_1+1,\ldots,k 
\end{equation}
Since $k$ is odd, (\ref{eq1.1}), (\ref{eq1.2}) and (\ref{eq1.3}) give a
contradiction. Note that the same contradiction holds if $k_2=0$.
Hence, if $k_1$ is odd $G$ cannot have an even $F$-orientation. Similarly,
(reversing the roles of (\ref{eq1.2}) and (\ref{eq1.3})) if $k_2$ is odd then
$G$ cannot have an odd $F$-orientation.
\end{proof}

\begin{cor}\label{claim2.1.2}
Let $G$ be a graph with a $1$-factor $F$ and an odd $F$-set.
Then $G$ cannot have both an odd $F$-orientation and an even $F$-orientation.
\end{cor}

\begin{proof}
In the notation of Lemma \ref{claim2.1.1}, since $k$ is odd either $k_1$ is odd or $k_2$ is odd.
Then the result follows directly from Lemma \ref{claim2.1.1}.
\end{proof}
%


\

The {\em Wagner graph $W$} is the cubic graph having vertex set
$V(W)=\{1,\ldots,8\}$ and edge set $E(W)$ consisting of the edges of the cycle $C=(1,\ldots,8)$ and the chords $\{(1,5),(2,6),(3,7),(4,8)\}$.

Let $C_1$ and $C_2$ be cycles of $G$ such that both include the pair of distinct independent edges $e=(u_1,u_2)$ and $f=(v_1,v_2)$. We say that $e$ and $f$ are {\em skew relative to $C_1$ and $C_2$} if the sequence $(u_1,u_2,v_1,v_2)$ occurs as a subsequence in exactly one of these cycles. Equivalently, we may write, without loss of generality, $C_1:=(u_1,u_2,\ldots,v_1,v_2,\ldots)$ and $C_2:=(u_1,u_2,\ldots,v_2,v_1,\ldots)$ i.e. if the cycles $C_1$ and $C_2$ are regarded as directed cycles, the orientation of the pair of edges $e$ and $f$ occur differently.

\begin{lemma}\label{claim2.6}\cite{ALS}
Let $F_1:=\{(1,5)$, $(2,6)$, $(3,7)$, $(4,8)\}$ and $F_2:=\{(1,2)$, $(3,4)$, $(5,6)$, $(7,8)\}$ be $1$-factors of the Wagner graph $W$. Set $e:=(1,8)$ and $f:=(4,5)$. Then the Wagner graph $W$ satisfies the following:

\begin{enumerate}[(i)]
  \item $W$ is $1$-extendable
  \item $W -\{e,f\}$ is bipartite and $1$-extendable (i.e. $W$ is near bipartite). 
  \item $W$ has an even $F_1$-orientation and an odd $F_1$-orientation.
  \item $W$ is Pfaffian.
  \item $W$ has no even $F_2$-orientation.
  \item There exist no pair of $F_1$-alternating cycles relative to which $e$ and $f$ are skew.
  \item The edges $e$ and $f$ are skew relative to the $F_2$-alternating cycles $C_1=(1,\ldots,8)$ and $C_2=(1,2,6,5,4,3,7,8)$.
\end{enumerate}
\end{lemma}


%
%
%
%
%
%
%
%
%
%
%

\begin{definition}$(${\em Generalized Wagner graphs $\cal{W}$}$)$\label{genwagner}
A graph $G$ is said to be a {\em generalized Wagner graph} if

\begin{enumerate}[(i)]

\item $G$ is $1$-extendable;

\item $G$ has a subset $R:=\{e,f\}$ of edges such that $G-R$ is 1-extendable and bipartite (i.e. $G$ is near bipartite);
\item $G-R$ has a $1$-factor $F$ and $F$-alternating cycles $C_1$ and $C_2$ relative to which $e$ and $f$ are skew.
\end{enumerate}
The set of such graphs is denoted by $\cal{W}$, and a $1$-factor $F$ of $G$ satisfying (iii) is said to be a $W$-factor of $G$.
\end{definition}

Note that, if we say that $G \in \cal{W}$, we will assume the notation of Definition \ref{genwagner} i.e. that $F$ is a $\cal{W}$-factor of $G$ and $R$, $C_1$ and $C_2$ are as described in Definition \ref{genwagner}(ii) and (iii), respectively.

Recently the authors proved in \cite{ALS}
the following result:

\begin{theorem}\cite{ALS}\label{lem1_6}
Let $G$ be a $1$-extendable graph containing a $1$-factor $F$
such that $G$ has no even $F$-orientation. Then $G$ contains an $F$-central subgraph
$G_0$ such that $G_0 \in \cal{W}$ and $F^*$ is a $\cal{W}$-factor of $G_0$.\qed
\end{theorem}

Note that in a companion paper \cite{ALS}, we complete this characterization
in the case of regular graphs, graphs
of connectivity at least four and of $k$-regular graphs for $k\ge3$.
Moreover, note that if $G_0 \in \cal{W}$ then $G_0$ is near bipartite.
Furthermore $F^*$ is the $1$-factor of $G_0$ induced by $F$ in the obvious way.

\section{Bad graphs}\label{bad}

In this section we introduce the definition of {\em bad} graphs and we study their relation with even and odd $F$-orientations. The results contained in this section will be fundamental in proving our main Theorem \ref{BadNonPfaf} 



\begin{definition}\label{badgraph2}$(${\em Bad Graph}$)$
A graph $G$ is said to be {\em bad} if $G$ contains a $1$-factor $F$ such that:
\begin{enumerate}[(i)]
\item $G$ has a zero-sum $F$-set $\cal A$;

\item $G$ has an orientation in which exactly an odd number of elements of $\cal A$ are evenly $F$-oriented (the other elements of $\cal A$ being oddly $F$-oriented).
\end{enumerate}
\end{definition}

This definition is equivalent to the one of \emph{intractable set of alternating circuits} given by Little and Rendl in \cite{LR}. We will prove, in 
Theorem \ref{BadNonPfaf} that a graph is bad if and only if it is non-Pfaffian, which corresponds to the equivalent result proved by Little in \cite{L73}, using linear algebra arguments.

\begin{definition}\label{badgraph1}$(${\em Simply Bad Graph}$)$
Let $G$ be a graph. $G$ is said to be \emph{simply bad} if $G$ contains a $1$-factor $F$ such that:
\begin{enumerate}[(i)]

\item $G$ has an odd $F$-set $\cal A$;

\item $G$ has an $F$-orientation in which each element of $\cal A$ is evenly $F$-oriented.
\end{enumerate}
\end{definition}


\begin{remark}\label{simplybadimpliesbad}
By definition, a simply bad graph is also bad. Definitions of bad and simply bad are, in fact, equivalent (this follows from Proposition \ref{thm1_9} and Theorem \ref{BadNonPfaf}).
\end{remark}

\begin{lemma}\label{claim1_4}
The graphs cubeplex $\Gamma_1$, twinplex $\Gamma_2$  and $K_{3,3}$ are simply bad.
\end{lemma}

\begin{proof}
$(i)$ {\em $\Gamma_1$ is simply bad:}

Let $F_1:=\{(a,d),(b,g),(i,c),(j,e),(h,k),(f,l)\}$. Note that $\Gamma_1 \setminus \{(c,i), (f,l)\}$ is bipartite and $\Gamma_1$ is not bipartite.
Hence $\Gamma_1$ is near bipartite. Let $\cal A$ be the set of $F_1$-alternating cycles defined by:

$$
\begin{array}{ll}
  C_1:=(a,d,c,i,j,e,f,l,k,h,g,b,a) & C_2:=(a,d,e,j,k,h,i,c,b,g,f,l,a) \\
  C_3:=(b,g,f,l,k,h,i,c,b) & C_4:=(a,d,c,i,j,e,f,l,a) \\
  C_5:=(a,d,e,j,k,h,g,b,a) &
\end{array}
$$

%

Thus, $\cal A$ is an odd $F_1$-set in which each element of $\cal A$ is evenly $F_1$-oriented.
Hence, $\Gamma_1$ is simply bad.

$(ii)$ {\em $\Gamma_2$ is simply bad:}

Note that $\Gamma_2$ may be obtained from the Petersen graph by subdividing two fixed edges at a
maximum distance apart and then joining the vertices of degree $2$ by an edge.
Notice that $\Gamma_2 \setminus \{(i,j), (e,f)\}$ is bipartite and hence $\Gamma_2$ is near bipartite.
Let $F_2:=\{(a,b),(c,d),(e,f),(g,h),(i,j),(k,l)\}$.
Let $\cal A$ be the set of $F_2$-alternating cycles defined by:

$$
\begin{array}{ll}
  C_1:=(a,b,f,e,l,k,g,h,d,c,j,i,a) & C_2:=(h,g,f,e,l,k,j,i,h) \\
  C_3:=(a,b,f,e,d,c,j,i,a) & C_4:=(a,b,c,d,h,g,k,l,a) \\
  C_5:=(a,b,c,d,e,f,g,h,i,j,k,l,a) &
\end{array}
$$

%
%
%
%

Thus, $\cal A$ is an odd $F_2$-set in which each element of $\cal A$ is evenly $F_2$-oriented.
Hence, $\Gamma_2$ is simply bad.

$(iii)$ {\em $K_{3,3}$ is simply bad:}

Finally, it is easily shown that $K_{3,3}$ is simply bad (see Figure \ref{fig2}).

\begin{figure}[h]
\begin{center}
\includegraphics[scale=0.25]{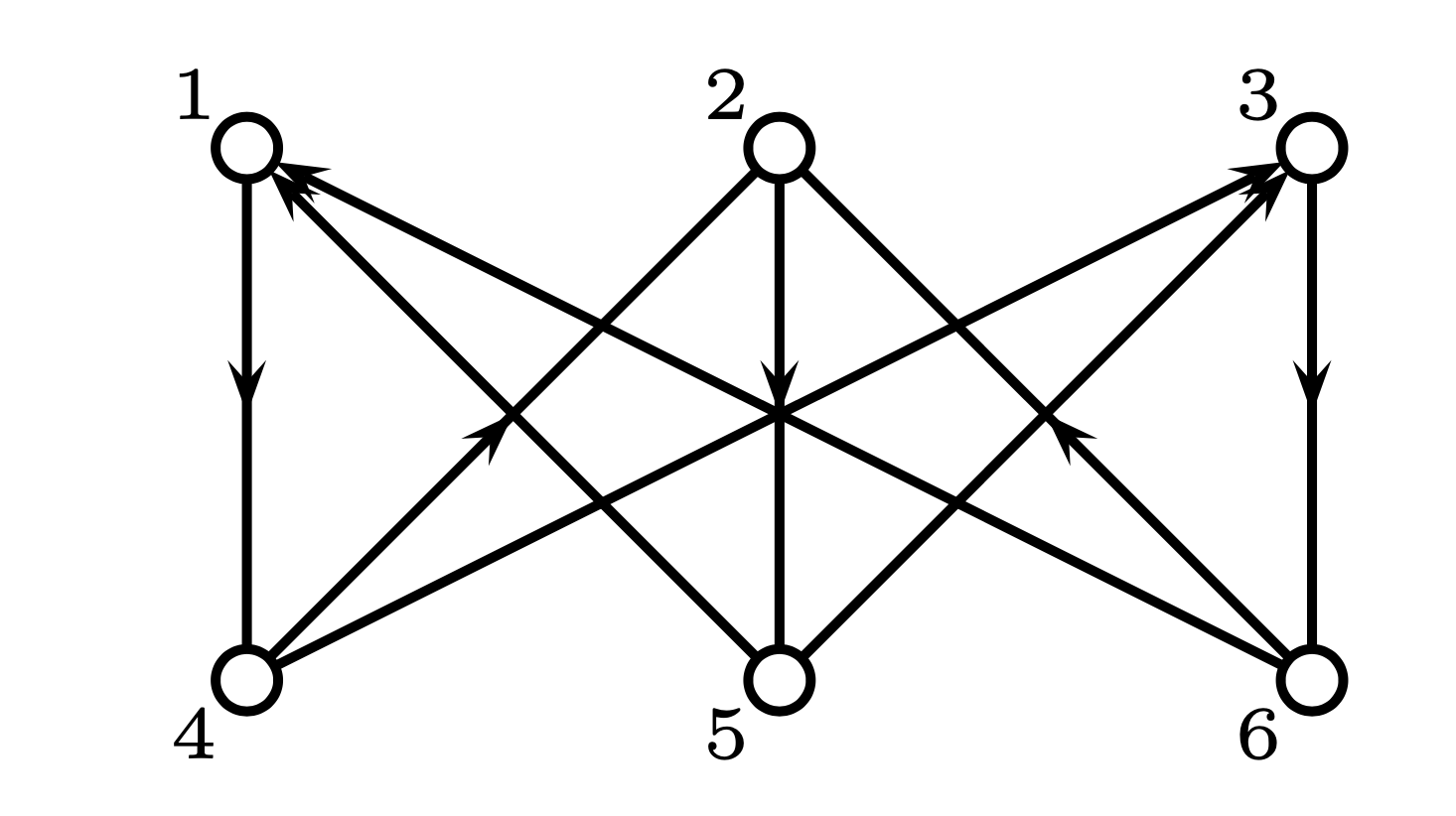}
\end{center}
  \caption{Orientation of the graph $K_{3,3}$ }\label{fig2}
\end{figure}

Using the notation of $(i)$ and $(ii)$, set $F_3:=\{(1,4), (2,5), (3,6)\}$
and ${\cal A}:= \{C_i | i=1,2, \ldots, 5\}$ where

$$
\begin{array}{ll}
  C_1:=(1,4,2,5,3,6,1) & C_2:=(1,4,3,6,2,5,1)\\
  C_3:=(1,4,2,5,1) & C_4:=(1,4,3,6,1) \\
  C_5:=(2,5,3,6,2) &
\end{array}
$$

%
%
%
%

The proof follows immediately.
\end{proof}


\

In the following lemmas we examine the relations between even subdivision, reducibility and simply bad graphs.

\begin{lemma}\label{evensubBad}
An even subdivision $H$ of a simply bad graph $G$ is also simply bad.
\end{lemma}

\begin{proof}
Let $F$ be a $1$-factor of $G$, let $\cal A$ be an odd $F$-set and $\overrightarrow{G}$ an orientation of $G$ in which
all elements of $\cal A$ are oddly oriented. Let $F^*$ be the $1$-factor of $H$ naturally induced by $F$,
in which from each path $P_e$ in $H$ which replaced an edge $e \in E(G)$, alternating edges are chosen into $F^*$
according to $e$ belonging to $F$ or not.
Similarly, $\cal A$ induces a set of cycles ${\cal A}^*$ in $H$, in which each edge of a cycle of $\cal A$ which
had been replaced by a path in $H$, is replaced by that same path in the corresponding cycle in ${\cal A}^*$.
Finally, $\overrightarrow{G}$ induces an orientation $\overrightarrow{G^*}$ in $H$ in which every path $P_e$ of $H$
which replaced an edge $e$ from $G$, has all edges oriented in correspondence to the orientation of $e \in \overrightarrow{G}$.
Since $H$ is an even subdivision, by definition ${\cal A}^*$ turns out to be an odd $F^*$-set and $\overrightarrow{G^*}$ turns out to be an orientation in which every cycle of ${\cal A}^*$ is evenly $F^*$-oriented, so $H$ is simply bad.
\end{proof}

%

%
%
%
%


%

\begin{definition}\label{orientfunc}
Let $\vec{G}$ be an orientation of $G$. We define a $(0,1)$-function $\omega:=\omega_{\vec{G}}$ on the set of paths and cycles of $G$ as follows:

(i) For any path $P:=P(u,v)=(u_0,\ldots,u_n)$
$$
\omega(P) \, := \, |\{i \, : \, [u_i,u_{i+1}] \in E(\vec{G}), 0\le i \le n-1 \}| \, (\mbox{mod} \, \, 2) \,.
$$

Note that $\omega(P(u,v))\equiv \omega(P(v,u)) + n \, (\mbox{mod} \, \, 2)$.

(ii) For any cycle $C=(u_1,\ldots,u_n,u_1)$
$$
\omega(C) \, := \, |\{i \, : \, [u_i,u_{i+1}] \in E(\vec{G}), 0\le i \le n-1 \}| \, (\mbox{mod} \, \, 2) \,;
$$
where the suffixes are integers taken modulo n.

We say that $\omega$ is {\em the orientation function associated with $\vec{G}$}.
\end{definition}

\begin{lemma}\label{claimA_1}
Suppose that $G$ is a graph which is simply reducible to a graph $H$. Then if $H$ is simply bad, $G$ is simply bad.
\end{lemma}

\begin{proof}
Suppose that $H$ is obtained from $G$ by contracting the cycle $D$ to a vertex $u$, where for some integer $k$ ($k \ge 1$)
$D:=(u_1,u_2, \ldots, u_{2k+1})$.

Suppose that $H$ is simply bad. Let $F$ be a $1$-factor of $H$ such that $H$ contains an odd $F$-set $\cal A$ and $H$ has an
$F$-orientation $\overrightarrow{H}$ in which each element of $\cal A$ is evenly oriented. Let $\omega$ be the associated orientation
function. Suppose that $e_i:=(u,v_i)$, $i=1,2, \ldots, 2k+1$ are a subset of the edges incident to $u$ such that $e^*_i:=(u_i,v_i)$ are edges in $G$ (we will assume that such edges exist and this makes no difference to the argument).
We may assume that $e_1 \in F$. Set $F_1:=\{(u_{2i},u_{2i+1} | i = 1,2, \ldots, k)\}$ and
$F_2:=F_1 \cup \{F \setminus \{e_1\} \cup \{e^*_1\} \}.$
Thus, $F_2$ is a $1$-factor of $G$. Now define an $F_2$-orientation $\overrightarrow{G}$ of $G$ with orientation function $\omega_2$, as follows:

\begin{enumerate}[(i)]
\item $\omega_2(a,b):=\omega(a,b)$ for each $(a,b) \in H \setminus \{u\}$;
\item $\omega_2(u_i,v_i):=\omega(u,v_i)$ for $i=1, 2, \ldots, 2k+1$; 
\item $\omega_2(u_{i},u_{i+1}):=1$ for $i=1, 2, \ldots, 2k+1$ (indices taken modulo $2k+1$).
\end{enumerate}

Let $C_j$ be a typical $F$-alternating cycle of $H$ containing $e_1$ and $e_j$. Then there is a natural one to one correspondence
with $F_2$-alternating cycles $C^*_j$ in $G$. Thus set $C^*_{2i}$ to be the $F_2$-alternating cycle in $G$ obtained from $C_{2i}$
on replacing the path $(v_1,u,v_{2i})$ by $(v_1,u_1,u_{2k+1},u_{2k},\ldots, u_{2i},v_{2i})$.
Similarly set $C^*_{2i+1}$ to be the $F_2$-alternating cycle obtained from $C_{2i+1}$ by replacing
$(v_1,u,v_{2i+1})$ by $(v_1,u_1,u_2,\ldots, u_{2i+1},v_{2i+1})$.
By definition $w(C^*_j)=w(C_j)=0$.

Let ${\cal A}^*$ be the set of $F_2$-alternating cycles which is obtained form $\cal A$ by replacing each $C_j$ by $C^*_j$.
Thus each element of ${\cal A}^*$ is evenly $F_2$-oriented in $\overrightarrow{G}$. Furthermore, consider the modulo $2$ sums
of the cycles in ${\cal A}^*$. Thus this is an Eulerian graph contained in $D$ (since $\cal A$ is an odd $F$-set) and hence
is a union of even cycles contained in $D$. Hence since $D$ is an odd length cycle, this Eulerian graph is empty.
Thus, ${\cal A}^*$ is and odd $F_2$-set. Hence, $G$ is simply bad.
\end{proof}


\begin{lemma}\label{claimA_2}
Let $G$ be a graph. If $G$ contains a simply bad subgraph $J$ such that $G-V(J)$ has a $1$-factor, then $G$ is simply bad.
\end{lemma}

\begin{proof}
Let $J$ be as in the statement.  Since $J$ is simply bad, $J$ has a $1$-factor $F$ such that $J$ contains
and odd $F$-set $\cal A$ and $J$ has an $F$-orientation in which each element of $\cal A$ is evenly oriented.
Now set $F_2:=F \cup F_1$ where $F_1$ is a $1$-factor of $G-V(J)$. Thus $G$ contains $\cal A$ and $\cal A$ is and odd $F_2$-set
and in $G$, $\cal A$ has the induced $F_2$-orientation in which each element of $\cal A$ is evenly $F_2$-oriented.
\end{proof}

\begin{lemma}\label{claim1_3}
If $G$ is reducible to $H$ and $H$ is simply bad then $G$ is simply bad.
\end{lemma}

\begin{proof}
	
 It is an immediate consequence of Lemmas \ref{claimA_1} and \ref{claimA_2}.
 \end{proof}


\

\section{Equivalence between Bad and Non-Pfaffian Graphs }\label{maintheorem}


In this section we prove our main characterization 
Theorem \ref{BadNonPfaf}. 
Firstly we need the following lemma which relates Pfaffian graphs to even $F$-orientations, and an accessory characterization
of non-Pfaffian graphs in terms of simply bad graphs (c.f. Proposition \ref{thm1_9}).

\begin{lemma}\label{lem1_8}
Let $G$ be a non-Pfaffian graph containing a $1$-factor $F$.
Suppose that $G$ has an even $F$-orientation. Then $G$ is simply bad.
\end{lemma}

\begin{proof} 
We use the proofs of Lemma 8.3.1 and Theorem 8.3.8 contained in \cite{LP}.

Let $G$ be a non-Pfaffian graph with a $1$-factor $F$ such that $G$ has an even
$F$-orientation $\overrightarrow{G}$.

By Theorem 8.3.7(4) in \cite{LP} there is a set of $1$-factors $F_1, F_2, \ldots, F_r$ ($r > 0$)
of $G$ such that

\begin{equation}\label{eq1}
\sum_{j=1}^{r} F_j \equiv 0 \, (\mbox{mod} \, \, 2)
\end{equation}

(i.e. each edge belongs to an even number of these $1$-factors)

and

\begin{equation}\label{eq2}
\sum_{j=1}^{r} \ell(F_j) \equiv 1 \, (\mbox{mod} \, \, 2)
\end{equation}

where for each $F_j$, $\ell(F_j)$ satisfies $sgn(F_j) = (-1)^{\ell(F_j)}$
and $\ell(F_j) \in \{0,1\}$. Here $sgn(F_j)$ denotes the sign of the $1$-factor
$F_j$.

Let $\cal A$ be the family of all $F$-alternating cycles formed from $F \Delta F_j$
for $j=1, 2, \ldots, r$ (where $\Delta$ stands for the symmetric difference).
Also let $k_j$ denote the number of $F$-alternating cycles
formed from $F \Delta F_j$. We may assume that the vertices of $G$ are labelled so that
$sgn(F)=1$. Hence, as in Lemma 8.3.1.:

\begin{equation}\label{eq3}
sgn(F) sgn(F_j) = sgn(F_j) = (-1)^{k_j}
\end{equation}

and thus, as in the proof of Lemma 8.3.8,
$\ell(F_j) \equiv k_j \mod 2$. Hence, from (\ref{eq2}),

\begin{equation}\label{eq4}
|{\cal A}| = \sum_{j=1}^r k_j = \sum_{j=1}^r \ell(F_j) \equiv 1 \, (\mbox{mod} \, \, 2).
\end{equation}

Furthermore consider the sum of the cycles in $\cal A$ modulo $2$.
If $e \notin F$ then, from (\ref{eq1}), $e$ is contained in an even number of
$F_j$ $(i=1,2, \ldots, r)$. Thus the modulo $2$ sum of the cycles in $\cal A$
is a subset of $F$. But since the modulo $2$ sum of cycles must be an Eulerian
graph, it follows that the modulo $2$ sum of cycles in $\cal A$ is zero.
Hence $\cal A$ is a simply bad $F$-set. Hence $G$ is simply bad.
\end{proof}

\

We give a characterization of non-Pfaffian graphs in terms of simply bad graphs and then use it to prove our main result which characterizes non-Pfaffian graphs in terms of bad graphs.


\begin{prop}\label{thm1_9}
Let $G$ be a graph. Then $G$ is simply bad if and only if it is non-Pfaffian.
\end{prop}

\begin{proof}
Let $G$ be a simply bad graph. From the definition of \emph{simply bad graph},
it follows that $k_2=k$ in Lemma \ref{claim2.1.1}, and $k$ is odd so $G$ has no odd $F$-orientation. Hence $G$ is non-Pfaffian.

Now suppose that $G$ is non-Pfaffian. There are two cases to consider: (i) $G$ has an even $F$-orientation where $F$ is a $1$-factor
of $G$; (ii) $G$ has no even $F$-orientation, for all $1$-factors $F$.

{\sc Case }(i). Let $G$ be a graph with an even $F$-orientation where $F$ is a $1$-factor
of $G$. Then, $G$ is simply bad by Lemma \ref{lem1_8}.

{\sc Case }(ii). Suppose that for all $1$-factors $F$, $G$ has no even
$F$-orientation.
Then, 
from Theorem \ref{lem1_6} and subsequent note,
for each $1$-factor $F$, the graph
$G$ contains an $F$-central subgraph $G_0$ which is near bipartite and non-Pfaffian.
Hence, from Theorem \ref{thm1_5}, $G_0$ contains a central
subgraph $J$ which is reducible to an even subdivision of $K_{3,3}$, $\Gamma_1$ or $\Gamma_2$.

By Lemma \ref{claim1_4}, $K_{3,3}$, $\Gamma_1$ and $\Gamma_2$ are simply bad,
and so is any even subdivision by Lemma \ref{evensubBad}.
Thus, applying Lemma \ref{claim1_3}, the subgraph $J$ is simply bad. Hence, applying Lemma \ref{claimA_2} twice,
both $G_0$ and $G$ are simply bad.
\end{proof}

\begin{theorem}\label{BadNonPfaf}
Let $G$ be a graph. Then $G$ is bad if and only if it is non-Pfaffian.
\end{theorem}

\begin{proof}
Let $G$ be a bad graph. From the definition of \emph{bad graph},
it follows that $k_2=k$ in Lemma \ref{claim2.1.1}, and $k$ is odd so $G$ has no odd $F$-orientation. Hence $G$ is non-Pfaffian.

Let $G$ be a non-Pfaffian graph. By Proposition \ref{thm1_9}, $G$ is simply bad and by Remark \ref{simplybadimpliesbad}, $G$ is bad.
\end{proof}

\

As mentioned earlier, this result was equivalently stated in Little and Rendl \cite{LR} and proved using linear algebra arguments in \cite{L73}. Here, we have given a graph theoretical proof and extended Fischer and Little's result \cite{FL} on near-bipartite graphs in terms of bad graphs.

\end{document}